\documentclass[11pt,table]{article}
\usepackage[utf8]{inputenc}
\usepackage{amsmath, amssymb}
\usepackage{algcompatible}
\usepackage{graphicx}
\usepackage{xcolor}
\usepackage{rotating}
\usepackage{tablefootnote}
\usepackage{xfrac}
\usepackage{fullpage}
\usepackage[noend]{algpseudocode}
\usepackage[linesnumbered,ruled,vlined]{algorithm2e}
\usepackage{subfigure}
\usepackage{caption}
\usepackage{multirow}
\usepackage{tabularx}
\usepackage{float}
\usepackage{comment}
\usepackage{booktabs}
\usepackage[toc,title,page]{appendix}
\usepackage[affil-it]{authblk} 
\usepackage{etoolbox}
\usepackage{lmodern}
\usepackage{url}

\usepackage{setspace} 
\usepackage[margin=1.0in]{geometry}

\definecolor{ao(english)}{rgb}{0.0, 0.5, 0.0}

\SetKwInput{KwInput}{Input}         
\SetKwInput{KwOutput}{Output}        
\SetKwInput{KwReturn}{Return}

\makeatletter
\patchcmd{\@maketitle}{\LARGE\@title}{\fontsize{16}{19.2}\selectfont\@title}{}{}
\makeatother

\title{Hyperconnected Megacity Parcel Logistics:\\
Parcel Routing and Containerized Consolidation}
\date{}
\author[1,3]{Sara Kaboudvand \thanks{Corresponding author}}
\author[1,2,3]{Benoit Montreuil}
\author[2,3]{Martin Savelsbergh}
\affil[1]{Physical Internet Center}
\affil[2]{Supply Chain and Logistics Institute}
\affil[3]{School of Industrial and Systems Engineering, Georgia Institute of Technology, Atlanta, GA}

\begin{document}
\maketitle
\begin{abstract}
    In high-speed hyperconnected parcel logistics, intermediate hubs play a critical role in reaching economies of scale through consolidating disperse flows of goods. However, resorting small-size parcels at every hub can increase hubs’ workload and parcels’ total travel time. Such resorting can be reduced by smart dynamic containerized consolidation, encapsulating together parcels sharing service level and a subsequent destination. In this study, we introduce an optimization model enabling to assess the impact of such consolidation in reducing the expected total pickup-to-delivery times over a parcel logistic network. We provide empirical results for a synthetic urban environment, contrasting consolidation performance over different network structures and demand patterns. Based on the results from the synthetic network experiments, we also assess the containerization potentials in a real-world logistic network instance. 
\end{abstract}
\section{Introduction} \label{sec:intro}

The rapid growth of e-commerce through the last decade has created new challenges for parcel logistics \cite{Savelsbergh2016}. On the one hand, more cost-effective solutions are needed for the companies to survive and thrive in this competitive service environment; On the other hand, customers' rising expectations narrow the pickup-to-delivery time window a logistic system can secure and, as such, the room for more efficient planning. One of the fundamental approaches to cost savings in the parcel logistics literature has been through economies of scale and flow consolidation.

The high-velocity flow of small-size goods from many origins to many destinations inherent to last-mile logistics encourages hub-based network structures for better consolidation and economies of scale. However, sorting a large number of parcels at (intermediate) hubs requires a significant investment in real estate, human, and machine resources.  Furthermore, in today's last-mile delivery systems, the actual transportation time is not the dominant component of the total in-transit time for many parcels. Instead, a significant portion of parcels' total in-transit time is related to the processing of parcels at hubs and, possibly, the waiting of parcels at hubs. From the economical perspective, long waiting and processing times at hubs not only impose unnecessary handling and storage costs but also may impede offering tight delivery services and thus limit the market share. 

If, on the other hand, parcels going in the same direction are smartly grouped into containers to bypass the sorting process at busy hubs, significant reductions in transit time will be achieved (as well as reductions in cost). Container consolidation refers to integrating disjoint parcels flow, all heading to some joint next destination, into larger volume shipments. Consolidating parcels into containers not only reduces the time and effort spent in material handling and sorting processes (by decreasing the number of parcel touches), it also reduces the chances of in-transit damages to the parcels. Importantly, containerization can free up sorting capacity at critical hubs as containers bypass the sorting process, and simplifies handling, loading, and unloading processes at the hubs. Containerization can bring numerous advantages to the cold chain industry as well. Multiple Compartment Vehicles (MCVs) are commonly used for shipping cold chain products, with each compartment providing a different temperature. Using containers (of potentially different sizes) instead of compartments will allow for more flexibility in assigning the available shipping space to different batches of goods with specific temperature requirements.

In this paper, we assess the impact of effectively routing and consolidating parcels into containers on lightening the sorting load at critical hubs and accelerating parcels' journey in the last-mile delivery system. Parcel routing and container consolidation are interconnected decision problems. Shipping parcels through shorter paths can reduce the total transportation time but increase the total in-transit time by increasing the overall handling time at the hubs (resulting from the lower chances for flow consolidation). As such, we study the joint parcel routing and container consolidation problems to assess its resulting savings in total transit and handling times in last-mile delivery systems.

Our numerical analysis shows that significant transit and handling time savings can be achieved through container consolidation, with up to 20\% in transit time savings and up to 80\% in handling time savings.  On the other hand, the analysis also shows that the savings from container consolidation strongly depend on (1) the network configuration, e.g., the number, capacity, and location of hubs and links between hubs, (2) the demand pattern, e.g., the number, size, and distribution of commodities, and (3) the ratio between the average commodity size and the container size. In this study, we consider a single container size; considering multiple container sizes may improve the benefits but is left for future research.

The remainder of the paper is organized as follows. In Section 2, we summarize previous studies related to the area of containerized consolidation. The joint parcel routing and container consolidation problem is introduced in section 3. Section 4 presents the integer programming formulation developed for modeling this problem, and section 5 summarizes the results of an extensive computational study. Finally, the paper is concluded in Section 7. %section 1
\section{Literature Review} \label{sec:lit}

Many studies in the parcel logistics literature have investigated the impact of economies of scale in saving delivery costs \cite{Rogerson1993}. One stream of research leveraging this factor relates to the Flow Consolidation problem (FCP) for which an extensive body of literature exists \cite{Min1990}. FCP is primarily studied in the Less than Truck Load (LTL) industry and seeks time-based or quantity-based aggregation of commodities flow into larger shipments at intermediate/inter-modal hubs to save on the shipping cost \cite{TYAN2003,BOOKBINDER2002}. Some FCP studies have investigated the trade-off between inventory holding and shipping cost considering freight consolidation \cite{Marklund2011,Mutlu2007}. FCP addresses the freight consolidation into trucks and does not consider shipping containers.

Maritime shippers were the very first users of standardized containers to facilitate inter-modal transportation \cite{pls2015}. Several studies in the maritime transportation literature propose optimization or simulation models to decide on the ideal number and size of containers for each maritime shipper aimed at minimizing containers empty spaces and the total carbon footprint \cite{Ming2014,Salam2016}. In the literature, this problem is referred to as the Container Loading Plan. Despite the dominant and growing use of containers in maritime transportation, it yet has some limitations. For example, not all ports are capable of handling containers, finding cargo for container backhaul is difficult, and some careers may not be able to take advantage of the gained speed in loading/unloading using containers \cite{Salam2016}.

Alonso is among those who studied the container loading plan within the trucking industry. He proposes a model that groups parcels into homogeneous pallets, which are then loaded into trucks. The goal of his model is to minimize the total number of trucks used \cite{Alonso2017}. The benefits of containerization are even more evident to the cold chain industry. Perishable goods usually require certain temperatures to preserve their qualities. In cold chain transportation, it is common to use Multi-Compartment Vehicles (MCVs), which have several isolated spaces that can accommodate goods with different temperature requirements. This practice is very similar to (and even can be improved by) the use of temperature-controlled containers of potentially different sizes that can carry goods with different perishability conservation. In 2018, Hartman outlined the basic structure of the MCV loading problem and proposed a multi-dimensional load planning model that embeds the items' cost allocation as well \cite{Hartman2018}.

The Freight Consolidation and Containerization Problem (FCCP), first introduced by Qin et al. in 2014 \cite{qin2014freight}, includes the shipments' routing into the containerization decision at the price of simplifying the container load planning. It seeks to assign shipments to routes as well as containers (potentially of varying sizes) so as to minimize total transportation cost \cite{Heeswijk2018, Nasiri2017}. The FCCP problem is sometimes confused with the Freight Consolidation Problem. The main difference between FCP and FCCP is that FCP is not concerned with containerized consolidation. In fact, FCP deals with consolidating shipping units to better fill the trucks, while FCCP addresses the problem of aggregating shipping units into containers with high fill rates so that trucks move with less empty space. In both problems, the main goal is to minimize the total vehicle miles traveled and the subsequent shipping cost and environmental harm.

The FCCP only considers a one-leg container consolidation to minimize the total transportation costs; it does not address a multi-leg consolidation where material handling cost and effort at intermediate hubs become specifically important. A thread of research can be found in the railway transportation management context, called Railroad Blocking Problem (RBP), which directly takes into account the advantage of bypassing intermediate railway terminals. RBP concerns routing of shipping cars and classifying cars with similar itinerary into blocks that can travel through several hubs together. The cars that get hooked together to shape a block may not have the same origin and destination but at least share parts of their path. Given this structure, cars and blocks in the blocking problem can be considered as equivalent to the parcels and containers in the parcel logistic problem.

Asad was among the first who formulated the railroad blocking problem. His proposed model aimed to minimize the total car and locomotive transportation cost and the yards classification cost. He considered a limit on the total flow passing through each blocking arc, but no restriction on the number of cars classified in each yard \cite{Asad1980}. In 1998, Newton et al. modeled railroad blocking as a network design problem by considering yards as nodes and blocks as arcs. Their model also took into account the yard capacity in terms of the maximum number of cars that can be classified at each yard and considered different priority classes for the cars. Besides the transportation and yard classification costs, Newton's model also considers the delay cost associated with the cars re-classification \cite{Newton1998}. Hassani et al. consider both users' and providers' cost and, as such, also include train-related constraints such as trains' weight and length and the maximum number of trains operating on each train-arc \cite{Hasany2013}. In a recent article in 2021, Yaghini et al. propose an arc-based formulation for the railroad blocking problem and explicitly model the assignment of blocks to the trains \cite{Yaghini2021}.

Despite its numerous similarities, there are several important differences between the railroad blocking problem and the containerized consolidation problem in parcel logistics. Unlike containers, there is no explicit limit considered on the blocks' capacity as they form simply by hooking several cars together. Moreover, in the railroad blocking problem, it is assumed that each block is shipped with only one train from its origin to its destination; As such, the blocks are not transferred between trains at intermediate yards. In contrast, in the containerized consolidation problem, containers might get crossdocked and transferred to a different mover at intermediate hubs. Because of this later difference, no crossdocking capacity is considered in the railroad blocking problem, while it needs to be addressed explicitly in the containerized consolidation problem.

To the best of the authors' knowledge, there is only one published study in the literature that addresses the containerized consolidation problem within the context of parcel logistics while considering its direct advantages on the handling effort. In 2014, Chayanupatkul and Hall published a technical report that presented their findings on this problem in a project for USC METRANS transportation center. They studied package routing for long-haul freight shipments while accounting for containerization to minimize transportation and sorting costs. They used a path-based formulation on a transformed network. In their transformed network, for any route, alternative routes are built, each of which bypassing the intermediate stations in some way. They solved the LP relaxation of their formulation as a multi-commodity network flow with no containerization constraint imposed and examined applying three different heuristic approaches for finding the IP solution \cite{chayanupatkul2004freight}

Our work is different from the study presented by Chayaupatkul et al. in several aspects. First, their study is focused on long-haul air transportation, which required abiding by certain flight schedules. This enforces specific modeling of parcel-to-container and container-to-flight assignment for their problem, increasing the problem complexity by creating solution symmetry. In last-mile logistics, however, there is no need (and many times no guarantee) to respect a fixed truck schedule. Therefore, without loss of generality and in favor of model simplicity, we consider a specific truck departure rate on each shipping arc per time unit. This allows for aggregating commodities flow and, as such, avoiding symmetry. Second, in our model, we assume commodities are indivisible, whereas in the problem addressed by Chayaupatkul et al., commodity flows are considered to be divisible, i.e., two units of the same commodity could be shipped on different paths. The divisibility of commodities makes the model more tractable by replacing the binary indicators with positive continuous variables. Third, we investigate, in fine detail, the relationship between network characteristics and customers' demand pattern on potential savings through containerized consolidation. This is quite important as there is never a single structure that works equivalently well for all logistics networks. Furthermore, one may find that the long-time savings associated with one setup may justify the costs associated with transitioning to different network configurations. Therefore, it is crucial to assess the benefits of the proposed model in several representative logistic setups to provide a benchmark to broader use cases.

This study proposes a new formulation for the parcel routing and containerized consolidation problem to assess the potential of containerized consolidation in saving origin-destination handling and transit time across a range of different network configurations considering different hub capabilities and demand characteristics. The following section provides more details on the problem we address and its underlying assumptions.  %section 2
\section{Problem Description} \label{sec:descr}

The joint parcel routing and containerized consolidation problem addressed in this study is to decide on the route that each commodity takes from its origin to its destination and where along its route it is consolidated with other commodities. The goal is to minimize the total transit time of the commodities, where transit time includes transportation time along the arcs and sorting and cross-docking time at the hubs. Commodities may need to deviate from their shortest path to get consolidated with other commodities into containers that can bypass the sorting process at some intermediate hubs. The selected route for commodities must satisfy their delivery service promise, the sorting and cross-docking capacities at the hubs, and the transportation capacities of the links connecting the hubs.

Each commodity is a tuple associated with an origin, destination, quantity, and delivery service, which are assumed to be known in advance. Commodities’ quantity refers to the number of parcels per hour that need to be shipped from commodity’s origin to commodity’s destination, and delivery service determines the time available for delivering a commodity at its destination. Commodities can be delivered at any time before their promised delivery due, and no penalties are incurred for early delivery. As such, a commodity can take any path that requires less or equal time compared to its time-to-delivery. We further assume that the shipments associated with a particular commodity are unsplittable, i.e., they should all follow the same path, while the items associated with a commodity can be put in different containers. 

\begin{figure}
\begin{center}
  \includegraphics[width=0.4\textwidth]{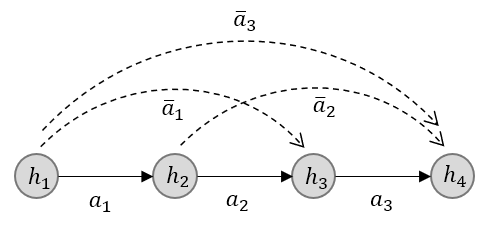}
  \caption{Physical and container arcs}
  \label{fig:modifiedNet}
\end{center}
\end{figure}

In order to incorporate containerization when making routing decisions, we use the network introduced in [1] but modify it to adapt to our problem's specific features. The modified network consists of a set of hubs and arcs. An arc can either be a physical arc representing a vehicle movement or a container arc, which represents a container movement. A container arc is defined by a sequence of hubs where the container is loaded at the hub at the tail of the arc and unloaded at the hub at the head of the arc and crossdocked at all intermediate hubs. Figure \ref{fig:modifiedNet} provides an example of a modified network with physical arcs $\{a_1, a_2, a_3\}$ shown with solid arrows and container arcs $\{\bar{a}_1,\bar{a}_2,\bar{a}_3\}$ are distinguished by dashed arrows.

\begin{figure}[H]
\begin{center}
  \includegraphics[width=0.4\textwidth]{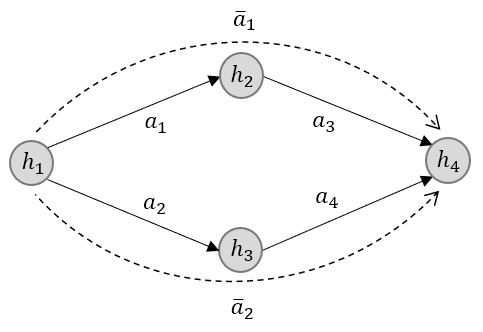}
  \caption{Different hub sequences corresponding to the same origin-destination sorting hubs}
  \label{fig:contaierArcs}
\end{center}
\end{figure}

Two container arcs may connect the same origin and destination sorting hubs, but have different hub sequences. For example, in Figure \ref{fig:contaierArcs}, the two container arcs $\bar{a}_1$ and $\bar{a}_2$ both connect sorting hubs $h_1$ and $h_4$, but correspond to different hub sequences which are $(h_1,h_2,h_4)$ and $(h_1,h_3,h_4)$ accordingly.

We assume that the physical arcs in the network are given. Furthermore, for each physical arc, we assume that the shipping capacity is provided in the form of the number of departures per unit time of a specific size vehicle from the hub at the tail of the arc. Given the vehicle departure frequency on a physical arc, i.e., the number of departures per unit time, we can also estimate the waiting time for a commodity at the hub at the tail of the arc. For example, with three vehicles departing per hour from hub $h_1$ to traverse the physical arc $a=(h_1,h_2)$ and each having capacity for 200 parcels, the shipping capacity along $a$ equals 600 parcels per hour, and the waiting time at $h_1$ is estimated to be 10 minutes (as the average time between vehicle departures is 20 minutes). The total estimated waiting time for a commodity from its origin to its destination is the sum of the estimated waiting times of the physical arcs along its path.

The following section explains our proposed solution methodology for solving the containerized consolidation problem under different logistic network configurations and demand patterns.

  %section 3
\section{Methodology} \label{sec:method}

This section first explains the mathematical formulation developed for solving the parcel routing and containerized consolidation problem. Next it covers the basis and methodologies for generating a variety of realistic problem instances, and finally expands on the IP generation including the container arcs creation algorithms as well as the logical restriction imposed on problem instances to heuristically decrease size of the problem.    

\subsection{Mathematical Formulation} \label{mathModel}
In this section, we introduce an integer programming formulation for the parcel routing and containerized consolidation problem using the arc and path structures introduced in Section \ref{sec:descr}.  We propose a path-based formulation that seeks to select a feasible path for each commodity such that the total transit time of the commodities is minimized (where the transit time includes transportation time, sorting time, cross-docking time, and waiting time). 

Consider a network $N = (H, A \cup \bar A)$ where $H$ represents the set of hubs, $A$ is the set of physical arcs and $\bar{A}$ refers to the set of container arcs.  Each hub $h\in H$ has a sorting time $c^S_h$ (minutes per parcel), a sorting capacity $l_h^S$ (number of parcels per hour), a cross-docking time $c^X_h$ (minutes per container), and a cross-docking capacity $l_h^X$ (number of containers per hour).  Each physical arc $a \in A$ has a transportation time $c^T_{a}$ (minutes), an estimated waiting time $c^W_{a}$ (minutes), and a vehicle departure frequency $d_{a}$ (number of vehicle departures per hour).
Let $K$ denote the set of commodities where for each commodity $k \in K$, $q_k$ show the commodity's quantity (number of parcels per hour) and $P(k)$ denote the set of feasible paths from commodity's origin to its destination. Also let $P = \cup_{k \in K} P(k)$ show all the paths generated in the target network. Let $A(p)$ denote the set of physical arcs in path $p \in P$, and let $H^X(p)$ and $H^S(p)$ denote the set of cross-docking and sorting hubs along path $p \in P$, respectively. Given these definitions, a path $p$ is feasible for commodity $k \in K$ if 

\[
T_p = \sum_{a \in A(p)} (c_{a}^T + c^W_a) + \sum_{h \in H^X(p)} c^X_h + \sum_{h \in H^S(p)} c^S_h
\]
is less than the allowed pickup-to-delivery time corresponding to the commodity's service level. 

For each physical arc $a \in A$, let $C(a) \subseteq \bar A$ denote the set of container arcs that use physical arc $a$. For each container arc $\bar a \in \bar A$, let $P(\bar a)$ denote the set of paths that include container arc $\bar a$ and let $H^X(\bar a)$ denote the set of cross-docking hubs along container arc $\bar a$. Finally, let $q$ denote the container capacity (number of parcels) and $Q$ denote the vehicle capacity (number of containers). We assume that all containers and all vehicles have the same size. We also assume that all parcels are uniform in size, and thus it suffices to specify container capacity in terms of the number of parcels, i.e., there is no need to consider the cumulative weight and volume of parcels. Lastly, commodity flows are assumed to be unsplittable, i.e., each commodity has to be shipped through a single path.

The IP formulation has two decision variables: (1) binary variable $x_k^p$ indicating if path $p$ is assigned to commodity $k$, and (2) the integer variable $y_{\bar a}$  indicating the number of containers moved along container arc $\bar a$ per hour. The IP formulation does not assign parcels to containers but instead assigns paths to commodities so as to encourage container consolidation. The assignment of parcels to containers, if needed, can be achieved by post-processing of the solution.

\begin{align}
    \min \quad &\sum_{k \in K} \sum_{p \in P(k)} q_k T_p x^{p}_{k} \\
    %##################################
    st. \quad & \sum_{\bar a \in \bar A \; : \; h\in H^X (\bar a)} y_{\bar a} \leq l_h^X && \forall h\in H\\
    %##################################
    &\sum_{k \in K} \sum_{p \in P(k) \; : \; h \in H^S(p)} q_k x^p_k \leq l_h^S &&\forall h\in H\\
    %##################################
    &\sum_{k \in K} \sum_{p \in P(k) \cap P(\bar a)} v_k x^p_k \leq q y_{\bar a} &&\forall \bar a\in \bar A\\
    %##################################
    &\sum_{\bar a \in C(a)} y_{\bar a} \leq Q d_{a} &&\forall a\in A\\ 
    %##################################
    &\sum_{p \in P(k)} x^p_k = 1 &&\forall k\in K\\
    %##################################
    &y_{\bar a} \in Z^+ && \forall \bar a \in \bar A\\
    %##################################
    &x_k^p \in \{0,1\} && \forall k\in K , \; \forall p \in P(k) 
\end{align}

The objective function minimizes the total transit time of the commodities (transportation time, sorting time, cross-docking time, and waiting time). Constraints (2) and (3) guarantee that the sorting capacity and cross-docking capacity at hubs are respected. Constraint (4) counts the number of containers on each container arc. Constraint (5) ensures that the shipping capacity along physical arcs in terms of the number of vehicles is respected. Finally, Constraint (6) ensure that each commodity is assigned to a single path. Constraints (7) and (8) specify the domains of the decision variables.

The next subsection discusses the methodologies developed to generate various instances with different network configurations and demand patterns.

\subsection{Instance Generation} \label{instGen}

To examine the benefits of containerized consolidation we generate synthetic last-mile delivery systems. These synthetic delivery system are characterized by three components: 1) the network configuration, i.e., the number and location of different types of hubs in the system and the transportation links that connect them, 2) the demand, i.e., the number, the size, and the geographic properties of the commodities, and the 3) capacity, i.e., the sorting and cross-docking capacity at the hubs and the transport capacity along the links connecting the hubs. In what follows, these three components are described in more details.

\subsubsection{Network Configuration}

\begin{figure}[htp]
  \includegraphics[width=\linewidth]{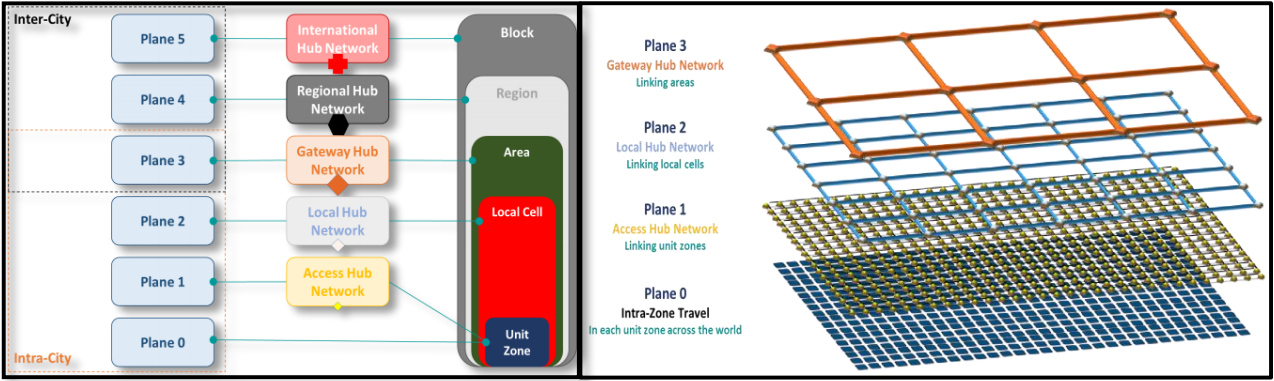}
  \caption{Multi-Tier Logistic Web [REFERENCE HERE]}
  \label{fig:multiTierWeb}
\end{figure}

We rely on the concept of a Multi-Tier Logistic Web introduced in \cite{montreuil2018urban} to generate networks with different configurations. Figure \ref{fig:multiTierWeb} illustrates the concept of a multi-tier logistic web. 

The multi-tier logistic web is built on top of a meshed network of \textit{unit zones} denoted as plane zero. A group of adjacent unit zones forms a \textit{local cell}, and a group of adjacent local cells forms an \textit{urban area}. At the higher tiers, a group of urban areas can form a \textit{region}, and a group of regions can form a \textit{country}. Except plain zero, each plane of the web corresponds to a set of hubs, with \textit{access hubs} located at plane one and serving the unit zones, \textit{local hubs} located at plane two and serving the local cells, and \textit{gateway hubs} located at plane three and serving the urban areas. There are also \textit{regional hubs} and \textit{international hubs} located at plains four and five and serve the regions and countries, respectively. 

\begin{figure}[htp]
\centering
  \includegraphics[width=0.8\textwidth]{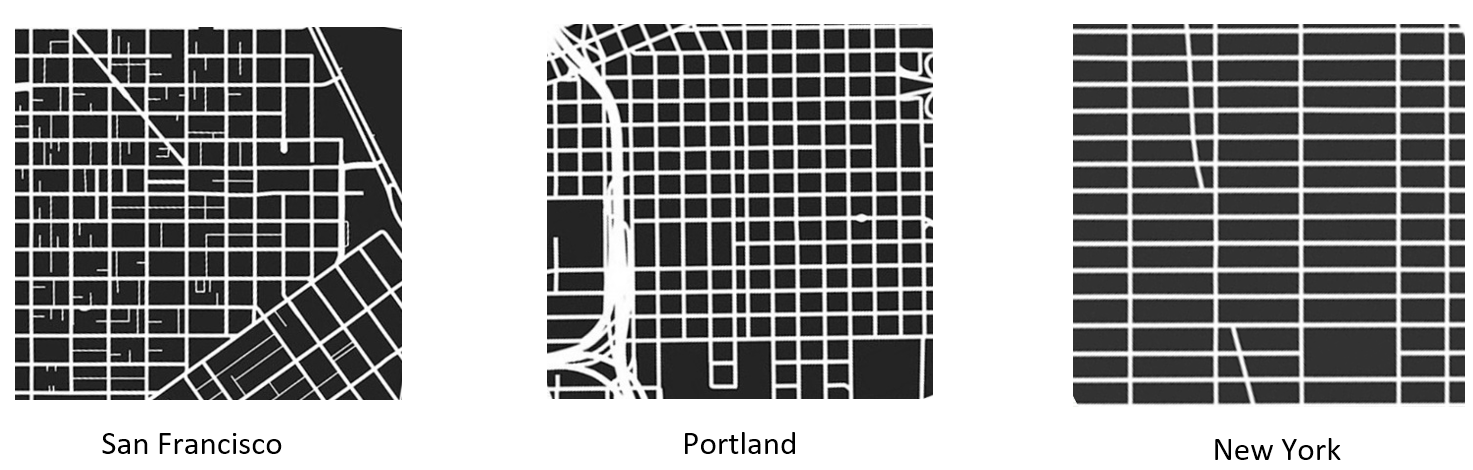}
  \caption{One Square Mile of the City Grids}
  \label{fig:cityGrids}
\end{figure}

\begin{figure}[htp]
\centering
  \includegraphics[width=1.0\textwidth]{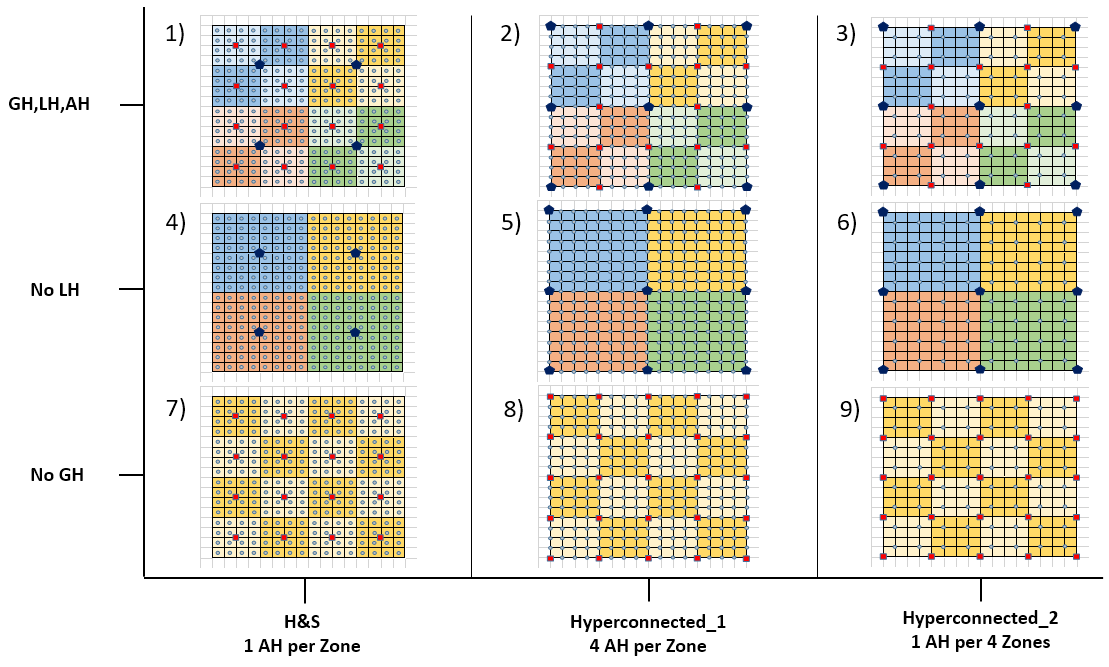}
  \caption{Experimental City Grids (COLORED)}
  \label{fig:exprCityGrids}
\end{figure}

Influenced by the layout of several US cities (see, for example, Figure \ref{fig:cityGrids}), we choose a simple grid structure for creating different network configurations. More specifically, we consider the nine grid structures depicted in Figure \ref{fig:exprCityGrids}. In each of these structures, the city is represented as a $16 \times 16$ grid. There are four urban areas, each indicated by a different color. Each urban area includes four local cells distinguished by lighter and darker themes from the same color, and each grid square represents a unit zone. Furthermore, Gateway, local, and access hubs are illustrated as blue pentagons, red squares, and gray circles, respectively.

Structure 1 represents a traditional hub-and-spoke network in which each unit zone is served by a single access hub located at the center of the zone, each local cell is served by a single local hub located at the center of the local cell, and each urban area is served by a single gateway hub located at the center of the urban area. Moreover, no direct shipment between neighboring access hubs and neighboring local hubs are allowed.

Structure 2 represents a hyperconnected logistic web \cite{montreuil2018urban} with each unit zone assigned to 4 access hubs located at the zone corners, each local cell served by 4 local hubs located at the cell corners, and each urban area served by 4 gateway hubs located at its corners. Moreover, unlike Structure 1, the hubs are not exclusively assigned to one area but are shared by adjacent areas. For example, each access hub is shared by up to four unit zones. Similarly, each local (gateway) hub is shared by up to four local cells (urban areas).  In this structure, direct shipment is allowed between adjacent access hubs in the same local cell and adjacent local hubs in the same area, which helps avoid unnecessary travel to higher tiers for shipping between nearby hubs.

Structure 3 is similar to structure 2 except for the number of access hubs and the assignment of unit zones to access hubs. Adjacent unit zones share access hubs, but each unit zone is only connected to one access hub. This structure will help to assess the impact of access hub density (number of access hubs) on containerized consolidation benefits for different demand patterns.

To assess the contribution of each tier of the network in facilitating containerized consolidation, structures 4, 5, and 6 are designed similar to structures 1, 2, and 3, respectively, but with the local hub tier removed. Moreover, structures 7, 8, and 9 mimic structures 1, 2, and 3, respectively, with the access hub tier removed.

\begin{figure}[htp]
\centering
  \includegraphics[width=0.7\textwidth]{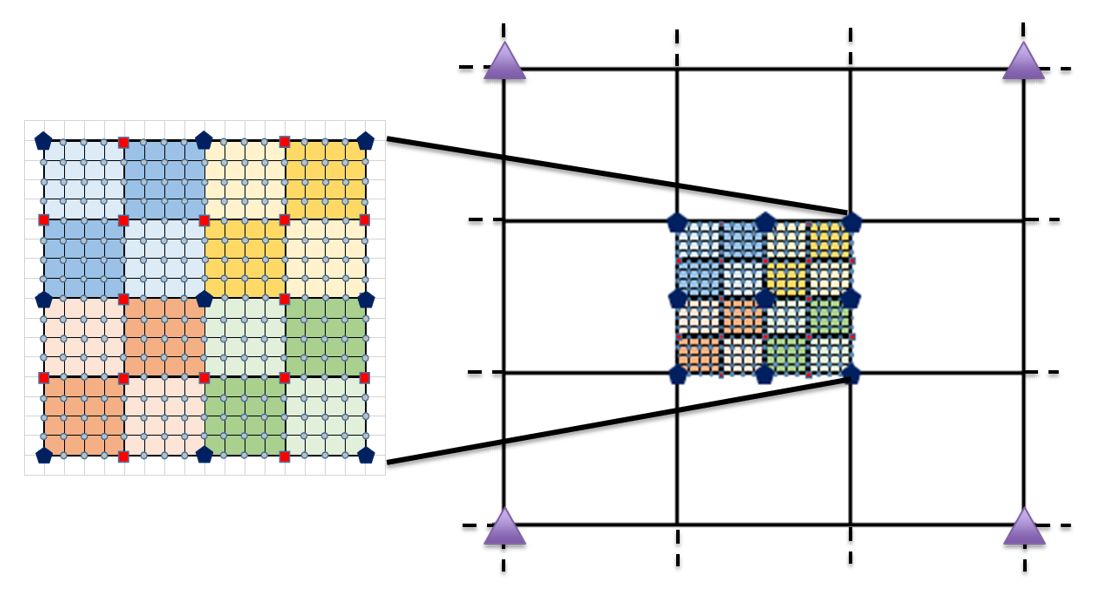}
  \caption{Embedding of a last-mile intracity delivery system in a larger intercity network (COLORED)}
  \label{fig:interIntraLink}
\end{figure}

To account for intercity parcel flows in a last-mile intracity delivery system, we also assume four regional hubs located at the northeast, northwest, southeast, and southwest corners of the city. Figure \ref{fig:interIntraLink} shows the embedding of a $16 \times 16$ last-mile intracity delivery system in a larger $48 \times 48$ grid structure where the purple triangles refer to the regional hubs. This embedding allows the modeling of different parcel flows in the last-mile delivery system. Commodities representing inbound parcels from other cities originate at one of the four regional hubs and are destined to one of the unit zones. Similarly, commodities representing outbound parcels to other cities originate from a unit zone and are destined to one of the regional hubs.

Lastly, the travel time along an arc between network nodes depends on the type of vehicle used on that arc, which may differ based on the type of node at the tail and the head of the arc as well as on the arc time distance. 

\subsubsection{Demand}

As mentioned earlier, the last-mile delivery system serves three types of demand:
\begin{itemize}
    \item \textit{Intracity parcels}: parcels originating in and destined for a unit zone within the city.
    \item \textit{Outbound Intercity Parcels}: parcels originating in a unit zone in the city, but destined for another city (which implies that they are destined for one of the regional hubs).
    \item \textit{Inbound Intercity Parcels}: parcels originating in another city and destined for a unit zone in the city (which implies they originate in one of the regional hubs).
\end{itemize}

Therefore, the unit zones and regional hubs are the potential pickup and delivery locations for commodities. To generate a realization of demand, we need to specify the number of commodities, number of parcels, and a pattern. A demand pattern defines how commodity origins and destinations are generated and thus determines the overall flow direction within a logistic network. 

\begin{figure}[htp]
\centering
  \includegraphics[width=0.7\textwidth]{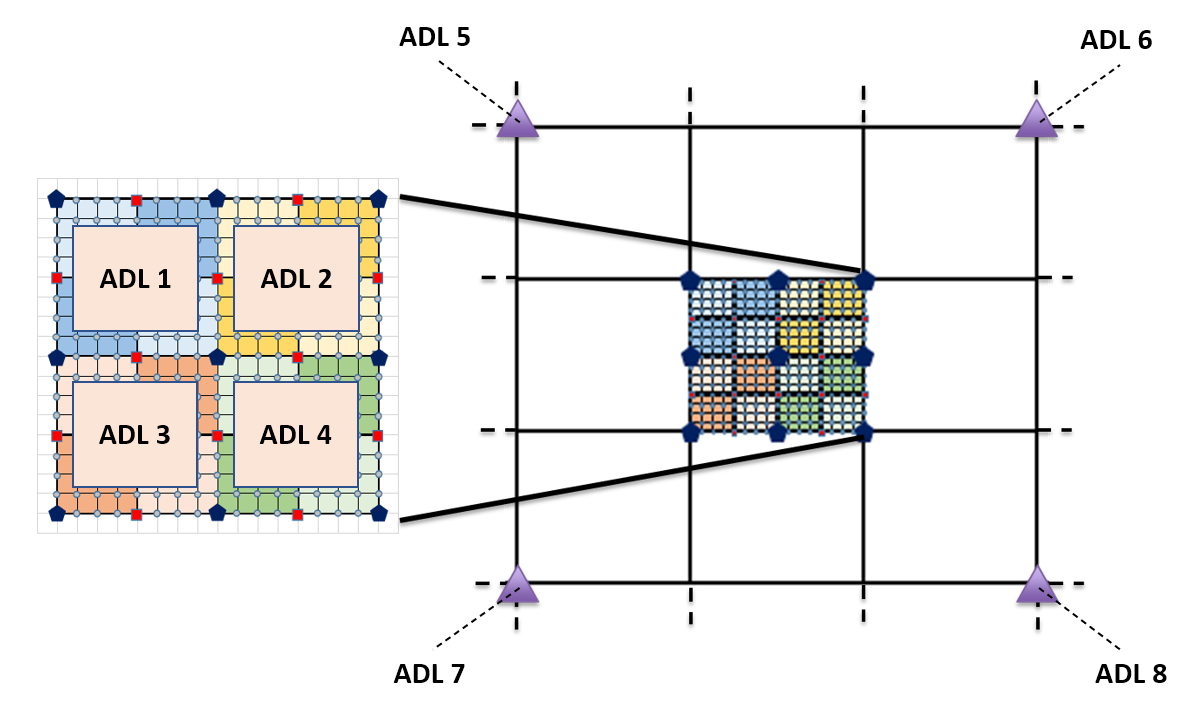}
  \caption{Clustering of demand points into Aggregate Demand Locations (ADLs) (COLORED)}
  \label{fig:adls}
\end{figure}

In order to define demand patterns, the potential pickup and delivery locations of commodities are grouped into Aggregate Demand Locations (ADLs). We consider one ADL for each urban area (a set of unit zones) and one ADL for each regional hub, summing up to the total of eight ADLs (see Figure \ref{fig:adls}). This classification is uniform across all network configurations.

The first step in defining a demand pattern is specifying the portion of the number of commodities over the three demand categories, i.e., intracity, intercity inbound, and intercity outbound. Each demand category implies eligible pickup and delivery ADLs. The list of eligible pickup and delivery ADLs for each demand category is summarized in Table \ref{tab:eligibleADLs}. For example, the four urban areas can serve as pickup ADLs for the intercity outbound demand category, and the four regional hubs can serve as delivery ADLs. After a portion of the total number of commodities was assigned to each demand category, for each category, we divide the pickup and delivery probabilities between its corresponding ADLs such that the sum of the pickup and the sum of the delivery probabilities equals 1. These pickup and delivery probabilities reflect demographic information about the ADL, e.g., relative population density and the relative ratio of business and residential occupancy.

\begin{table}[htp]
\small
\centering
\begin{tabular}{l|l|l}
\toprule
demand category & eligible pickup ADLs & eligible delivery ADLs\\
\midrule
Intracity & ADL1, ADL2, ADL3, ADL4 & ADL1, ADL2, ADL3, ADL4 \\ 
Intercity Inbound & ADL5, ADL6, ADL7, ADL8 & ADL1, ADL2, ADL3, ADL4 \\ 
Intercity Outbound & ADL1, ADL2, ADL3, ADL4 & ADL5, ADL6, ADL7, ADL8 \\ 
\bottomrule
\end{tabular}
\caption{Intracity ADLs' pickup and delivery probability for different demand patterns}
\label{tab:eligibleADLs}
\end{table} 

Given the number of commodities $n$, a demand category $c$ with associated fraction $f_c$, an associated ADL $i$ with pickup probability $p^o_{ci}$, and an associated ADL $j$ with delivery probability $p^d_{cj}$, the number of commodities of category $c$ with an origin in ADL $i$ and a destination in ADL $j$ is $\lceil n \cdot f_c \cdot p^o_{ci} \cdot p^d_{cj} \rceil$. 
Next, for each commodity with an origin in ADL $i$ and a destination in ADL $j$, we randomly select a location in ADL $i$ and a location in ADL $j$ (if $i = j$, we ensure that different locations are selected). After selecting the pickup and delivery location for the commodity, its size, i.e., number of parcels, is drawn from a triangular distribution with parameters $(a=1, c=m, b=2m)$, where $a$, $b$, and $c$ correspond to the minimum, maximum, and the average value, respectively, and $m$ represents the average number of parcels per commodity (i.e., the number of parcels divided by the number of commodities).

\begin{figure}[ht]
\centering
\subfigure[Uniform Demand Pattern]{
    \includegraphics[width=0.4\linewidth]{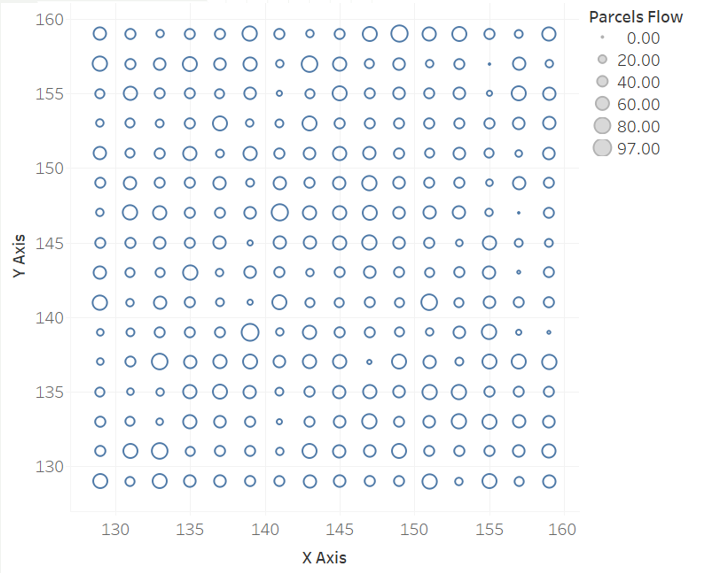}
    \label{fig:demandSubfig1}
}
\subfigure[Centric Demand Pattern]{
    \includegraphics[width=0.4\linewidth]{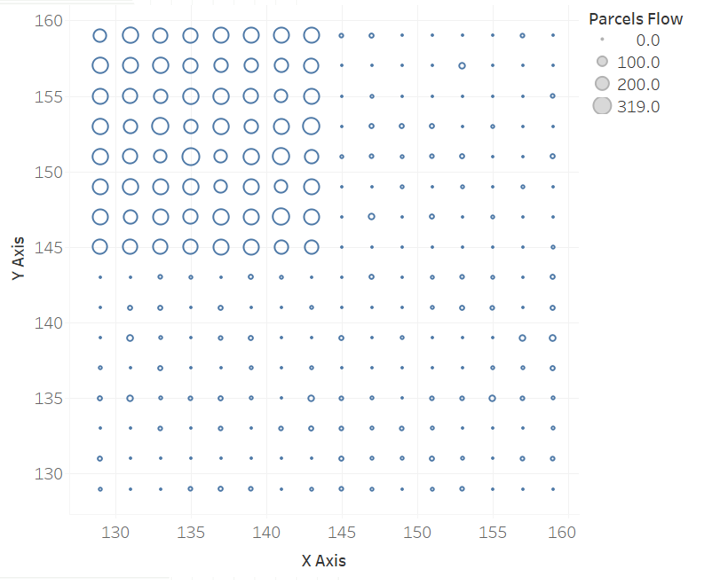}
    \label{fig:demandSubfig2}
}
\subfigure[Bi-polar Demand Pattern]{
    \includegraphics[width=0.4\linewidth]{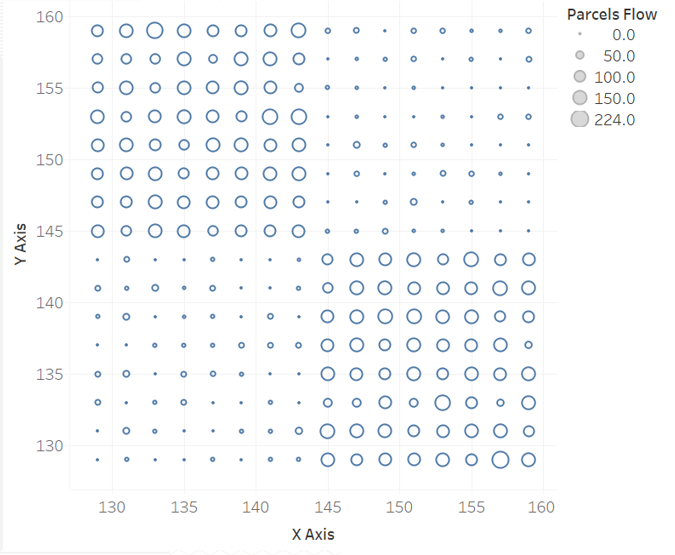}
    \label{fig:demandSubfig3}
}
\caption[Optional caption for list of figures]{Pickup and delivery density at different unit zones for different demand patterns}
\label{fig:flowPattern}
\end{figure}

We consider three demand patterns: (1) a uniform pattern: the pickup and delivery probabilities are distributed uniformly across the ADLs, (2) a centric pattern: 50\% of the pickups and 50\% of the deliveries occur in a single urban area, and (3) a bi-polar pattern: 50\% of the pickups originate in one urban area, and 50\% of the deliveries are destined for a different urban area. Figure \ref{fig:flowPattern} shows pickup and delivery density at different unit zones for different demand patterns.

After commodities are generated, their service promises are determined. The set of service promises to be considered, and the demand share for each service are decided in advance and given to the model as input. We compute a minimum doable delivery time for each commodity given its shortest path within the specific network configuration. Based on the minimum doable delivery time, we can then identify the set of commodities eligible for each service promise. Next, starting from the tightest delivery service promise and considering its target market share, we randomly select from its pool of eligible commodities to be assigned to that service. Extra commodities in each pool are considered for the next tightest service category as they are, by definition, included in the pool of eligible commodities for that service too. This approach guarantees that all commodities are assigned a delivery service while the target market share for each delivery service promise is met as best as possible. 

\subsubsection{Capacity \label{capacity}}

The last step in generating the synthetic last-mile delivery instances is to determine the shipping capacity along the physical arcs and the sorting and cross-docking capacity at the hubs. To do this, we first solve a Minimum Cost Multi-Commodity Network Flow (MCMCNF) problem to get an ``outline'' of how the commodities flow through the network. Since in real-world settings, not all commodities can and will be routed along the shortest path from origin to destination, we encourage a commodity to use multiple paths by assessing a penalty when the flow of a commodity on an arc is more than $\gamma$ of its size ($0 < \gamma < 1$). There is, however, no limit considered on the total flow passing through a physical arc. To find the worst-case minimum sorting capacity requirement, we assume that commodities are sorted at every terminal they visit along their path from origin to destination. Based on the industry's average, we consider an initial limit on the sorting capacity at different types of hubs. To ensure a feasible solution exists, we allow exceeding these limits, but the extra volume is penalized.

The following linear program represents the MCMCNF that we use for mapping commodity flows through the network. In this formulation, $N$ refers to the set of nodes, and $A$ refers to the set of physical arcs. Since commodities can also originate at or be destined for unit zones, the network nodes include both unit zones and hubs. We denote the set of nodes representing hubs as $N^H \subseteq N$. The set of commodities is denoted by $K$, and $o_k$, $d_k$, and $q_k$ represent the origin, the destination, and the size of commodity $k \in K$, respectively. Let $l_n^s$ denote the sorting limit at hub $n \in N^H$ and $c_a^T$ denote the travel time along physical arc $a \in A$. Decision variables $x_k^{a}$ indicate the flow of commodity $k \in K$ on physical arc $a \in A$. Variables $u_k^a$ represent the excess flow of commodity $k \in K$ on physical arc $a \in A$ (considering the $\gamma$ threshold) and variables $v_n$ represent the excess flow through hub $n \in N_H$. Finally, $M$ is a large value representing the unit-penalty imposed on the excess flow through the hubs and physical arcs.

\begin{align}
    \min \quad &\sum_{k \in K} \sum_{a \in A} (c_{a}^T x^{a}_{k} + M u_k^{a})) + \sum_{n \in N^H} M v_n \\
    %##################################
    st. \quad &\sum_{\substack{a \in A: \\a.tail=n}} x^a_k - \sum_{\substack{a \in A:\\a.head=n}} x^a_k = \begin{cases} q_k \quad\text{if} \quad n=o_k\\ -q_k \quad\text{if} \quad n=d_k \\ 0 \quad ow.\end{cases} &&\forall k,n\\
     %##################################
    &\sum_{\substack{a \in A:\\a.head=n}} \sum_{k \in K} x^a_k \leq l_n^s + v_n &&\forall n\in N^H\\
    %##################################
    & x^a_k \leq \gamma q_k + u^a_k && \forall k, a\\
    %##################################
    &x^a_k, u^a_k \in \mathbb{R}_{+} && \forall k,a \\
    &v_n \in \mathbb{R}_{+} && \forall n \in N_H
\end{align}

The objective function minimizes the penalty and the total shipping time in the system. Constraints (10) ensure flow conservation. Constraints (11) limit the flow of commodities through the hubs, and Constraints (12) capture the flow of commodities through the physical arcs. Constraints (13) and (14) specify the domains of the decision variables.

When the flow for each commodity $k$ across a physical arc $a$ is determined ($x^a_k$), the total flow along arc $a$ is computed as: $x_a = \sum_k x^a_k$. Next, considering a factor $\rho > 1$, the sorting capacity (in number of parcels) and the crossdocking capacity (in number of containers) at hub $n$, denoted by $l_n^s$ and $l_n^x$, respectively, are computed as $l_n^s = \rho \times v_n$ and $l_n^x = \lceil(\delta.l_n^s)/q \rceil$ where $\delta>0$ is a predefined constant and $q$ is the capacity of a container (in number of parcels). 

The shipping capacity along arc $a$ (in terms of the number of vehicles' departure from the tail of $a$ per hour) is computed as $\lceil(\rho\times x_a)/Q\rceil$, where $Q$ represents the uniform vehicles' capacity. In each scenario, depending on the network configuration and the demand pattern, a hub or arc might be assigned a zero capacity by the MCMCNF model, which implies that they are effectively removed from the network. 

The shipping time along a physical arc $a$ is computed considering two components: (1) transportation time, which depends on the arc distance, and (2) expected waiting for shipment at the tail of the arc, which is impacted by the departure frequency along that arc. For example, for an arc having two vehicles departing per hour, the trucks' inter-arrival time at its tail is 30 minutes on average, resulting in an estimated waiting time of 15 minutes for shipping.

Since, in the end, we desire to assign one path to each commodity, a final capacity adjustment is conducted to resolve any remaining infeasibilities in terms of hubs and arcs capacity. This final adjustment is made by solving the math model introduced in section \ref{mathModel} considering no consolidation and allowing but penalizing any excess flow through the physical arcs and the hubs. A sufficiently large penalty factor will assure that excess capacity at the hubs and along the physical arcs is used only if the model is infeasible otherwise. Such excess capacity is added to the base capacity, and the resulted feasible capacitated model is solved to derive numerical results. In the following subsection, we present the underlying assumptions for building IP instances that are fed into the mathematical model developed for parcel routing and containerized consolidation. 

\subsection{IP Generation}

As the number of feasible commodity paths becomes prohibitively large for even medium-size instances, we solve the IP formulation heuristically.  More specifically, we restrict the physical paths generated for a commodity in two ways: 
\begin{itemize}
    \item the length of a path cannot deviate more than a pre-specified factor from the length of the shortest path; and
    \item a path cannot contain more than a pre-specified number of intermediary hubs.
\end{itemize}

To feed the IP model, for each physical path, i.e., hub sequence, we need to generate a set of associated container paths, i.e., sorting hub sequences. In the pickup and delivery problems, the flow naturally becomes more disperse at the lowest tier of the network close to the pickup and delivery points. As such, except for relatively large-volume commodities, there are little chances for a commodity to get crossdocked through all intermediate terminals along its path. Because of that, to further reduce the size of the problem, we use Algorithm \ref{cArcs_alg} for generating container paths where we limit the number of terminals in a container arc to at most $K$.

\vspace{6pt} 
\begin{algorithm}[H]
\SetAlgoLined
\DontPrintSemicolon
\KwInput{A physical path $p = (n_0, n_1, n_2, \ldots, n_S)$ and a limit $K$ on the number of terminals in a container arc} 
\KwOutput{A set of container paths and a set of container arcs} 
$CP \leftarrow \emptyset$ \;
$CA \leftarrow \emptyset$ \;
$c \leftarrow (n_0)$ \;
$Expand(CP,CA,c,0)$ \;
return $CP$, $CA$ 
\caption{Container Path Generation}
\label{cArcs_alg}
\end{algorithm}

\begin{algorithm}[H]
\SetAlgoLined
\DontPrintSemicolon
\KwInput{ \; $CP$: collection of container paths \; $CA$: collection of container arcs \; $c$: partial container path \; $i$: index $i$ \; }
\tcp{$append(c,n)$: adds terminal $n$ at the end of partial container path $c$} 
\For{$j=i+1,\ldots, \min \{i+K+1, S\}$}{
    $c' \leftarrow append(c,n_j)$ \tcp*{expand partial container path} 
    $CA \leftarrow CA \cup \{(n_i, n_{i+1}, \ldots, n_j)\}$ \tcp*{add container arc}
    \uIf{$j = S$}{
        $CP \leftarrow CP \cup \{c'\}$ \tcp*{add container path}
    }
    \Else{
        $Expand(CP, CA, c', j)$ 
    }
}
\caption{$Expand(CP,CA,c,i)$}
\end{algorithm}

\vspace{12pt}

  %section 4
\section{Computational Study \label{sec:comp}}

This section aims to assess the potential benefits of containerized consolidation for reducing commodities' in-transit times in different network configurations and under a variety of demand patterns. For doing so, relying on the methodologies presented in Chapter \ref{sec:method}, we generate several synthetic last-mile delivery systems that resemble real-world urban logistic networks, but that also have characteristics that facilitate detailed analysis of the benefits of containerized consolidation. 
 
 Within the synthetic networks, we assume that each grid cell is a 2km$\times$2km unit zone. Furthermore, distances along the physical arcs are computed based on rectilinear motion, and transportation times are calculated based on pre-defined assumptions on vehicles' speeds through different types of arcs (shown in Table \ref{tab:movers'Speed}). For example, on an 18km arc connecting an access hub and a local hub, the vehicle speed is considered to be 30 km/hr resulting in a 36 minutes transportation time. We assume that 50\% of the commodities have a 5-hour delivery promise, and the remaining 50\% are due in 10 hours.

\begin{table}[!ht]
\scriptsize
\centering
\begin{tabular}{ p{5cm}|p{1.5cm}|p{1.5cm}|p{1.5cm}|p{2cm}}
\toprule
\multirow{2}{4em}{Arc (tail:head)} & \multicolumn{3}{c|}{Speed (km/hr) for distance} & \multirow{2}{4em}{Capacity (\#parcels)}\\
& 0-10 km & 10-20 km &  $>$20 km & \\
\midrule
UZ:AH/UZ:LH/UZ:GH/UZ:RH & 12 & 12 & 12 & 60 \\ 
\midrule
AH:AH/AH:LH/AH:GH/AH:RH & 20 & 30 & 45 & 300 \\ 
\midrule
LH:LH/LH:GH/LH:RH & 30 & 40 & 55 & 1000 \\ 
\midrule
GH:GH/GH:RH & 50 & 60 & 65 & 3500 \\ 
\midrule
RH:RH & 70 & 80 & 100 & 3500 \\ 
\bottomrule
\end{tabular}
\caption{Movers' Speed Assumption based on Type and Distance of Physical Arcs}
\label{tab:movers'Speed}
\end{table} 

All numerical experiments are performed with AMD EPYC processor (with IBPB) 2.50 GHz (2 processors), with 120 GB assigned RAM and Windows Server 2012 R2 standards, 64-bit operating system, x64-based processor. We set the optimality gap for all scenarios to \%0.01. 
 
\subsection{Analysis on Different Operational Setups} 

Our first set of computation experiments focuses on assessing the savings of containerized consolidation on transit time, i.e., pickup to delivery time, in a hyperconnected logistic web (Structure 2 -- see Figure \ref{fig:exprCityGrids}). We consider 10,000 individual parcels divided between 1,000 commodities and a uniform demand pattern across the city territory. Moreover, The estimated arc and node flows are multiplied by 1.3 to obtain arc and node capacities. The cross-docking capacity at a hub is assumed to be four times the sorting capacity divided by the container size. Furthermore, we assume the sorting process at a hub takes four times more time than the cross-docking process. Finally, we assume that each container can accommodate up to 40 parcels -- four times the average demand per commodity. This allows us to determine the truck capacities, in terms of the number of containers, on different types of arcs (i.e., the round down of the truck capacity in number of parcels divided by the container capacity in number of parcels).
When generating commodities' path, we allow up to 7 intermediate hubs and up to 5\% deviation from the shortest path. We also limit the number of alternative paths per commodity to at most 20. 

The results demonstrate that containerized consolidation can bring significant benefits in terms of transit time and handling time with savings equivalent to \%19.51 and \%71.21, respectively. The commensurate savings in handling effort also points to likely reductions in handling costs (e.g., due to a reduction in the required workforce). Furthermore, the model run time is quite promising (174.77 sec for reaching \%0.01 optimality), indicating the scalability of the proposed IP model. 
In what follows, we expand on a number of sensitivity analysis over different operational characteristics to provide more insights on the containerized consolidation benefits.

\begin{figure}
\begin{center}
  \includegraphics[width=\textwidth]{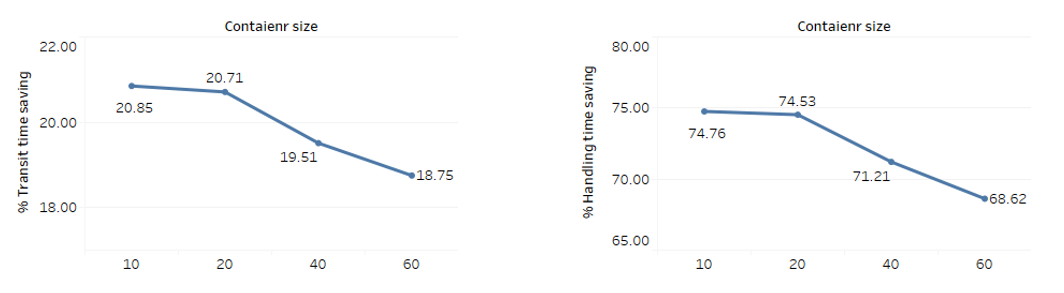}
  \caption{Sensitivity analysis over container size}
  \label{fig:contSize}
\end{center}
\end{figure}

\begin{figure}

\begin{center}
  \includegraphics[width=0.8\textwidth]{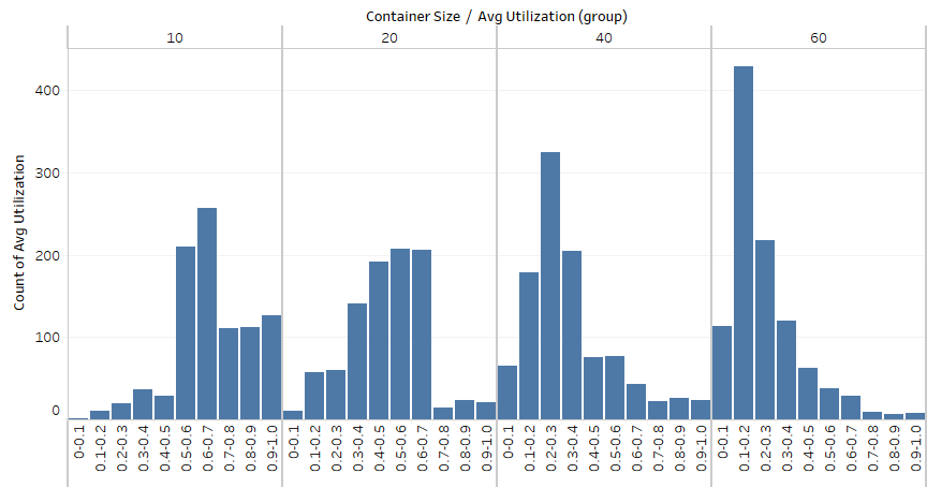}
  \caption{Utilization for different container sizes}
  \label{fig:contSizeUtil}
\end{center}
\end{figure}

Figure \ref{fig:contSize} shows the impact of container size on the transit and handling time savings achieved through containerized consolidation. As we see in the figure, there is no significant difference in the savings from containerization between container sizes 10 and 20. However, for larger container sizes, the savings from containerization decreases. This is because using larger containers for low volume flows will require more frequent container sorting. In other words, with smaller containers, there is more potential for consolidating commodities that can travel together over longer distances. Moreover, as the container size increases, the average containers utilization decreases (Figures \ref{fig:contSizeUtil}). Let $f_c$ represent the total flow shipped along container arc $c$ and let $q$ indicate container size in number of parcels. The average container utilization along container arc $c$ is then computed as follows:
\begin{align*}
    \text{avg. utilization along $c$} = 1- (\frac{\lceil \frac{f_c}{q}\rceil-\frac{f_c}{q}}{\lceil\frac{f_c}{q}\rceil}).
\end{align*}

\begin{figure}[htp]
\begin{center}
  \includegraphics[width=0.9\textwidth]{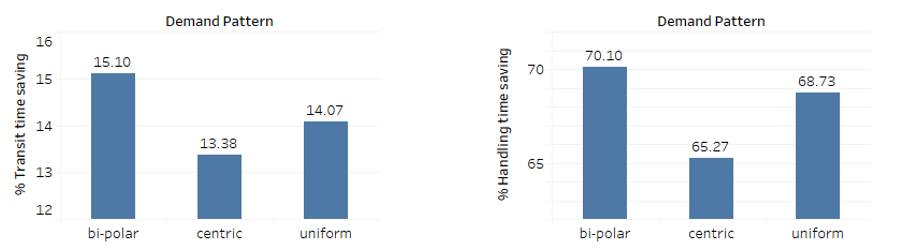}
  \caption{Percentage of transit and handling time savings for different demand patterns}
  \label{fig:demandPattern}
\end{center}
\end{figure}

Figure \ref{fig:demandPattern} shows the impact of the demand pattern on the savings from containerization. To make the comparison more revealing, we assume that the fraction of commodities representing intracity packages is \%80 -- rather than 50\% as in the base scenario -- and the remaining 20\% is equally divided between intercity inbound and intercity outbound commodities.  Table \ref{tab:demandPattern_Prob} shows the pickup and delivery probability distributions for the intracity ADLs in each demand pattern.  This implies, for example, that the fraction of commodities shipped from ADL 1 to itself in a centric demand pattern is equal to $0.8 \times 0.79 \times 0.79$, which is approximately 0.50. Similarly, when the demand pattern is bi-polar, approximately 50\% of the commodities have an origin in ADL 1 and a destination in ADL 4.  For intercity inbound and intercity outbound demand, we assume that commodities are uniformly distributed across the pickup and delivery ADLs.

\begin{table}[htp]
\scriptsize
\centering
\begin{tabular}{l|cc|cc|cc}
\toprule
\multirow{2}{2.5em}{ADL} & \multicolumn{2}{c|}{uniform demand} & \multicolumn{2}{c|}{centric demand} & \multicolumn{2}{c}{bi-polar demand}\\
& pickup pr. & delivery pr. & pickup pr. & delivery pr. & pickup pr. & delivery pr. \\
\midrule
ADL 1 & 0.25 & 0.25 & 0.79 & 0.79 & 0.79 & 0.07\\ 
ADL 2 & 0.25 & 0.25 & 0.07 & 0.07 & 0.07 & 0.07\\ 
ADL 3 & 0.25 & 0.25 & 0.07 & 0.07 & 0.07 & 0.07\\ 
ADL 4 & 0.25 & 0.25 & 0.07 & 0.07 & 0.07 & 0.79\\ 
\bottomrule
\end{tabular}
\caption{Intracity ADLs' pickup and delivery probability for different demand patterns}
\label{tab:demandPattern_Prob}
\end{table} 

Numerical results suggest that a bi-polar demand pattern results in the largest handling and transit time savings due to the natural flow convergence in this structure imposed by the overall flow direction. Furthermore, we see that a centric demand pattern results in the smallest transit time savings. The reason is that in this setup, the majority of commodities have a relatively short distance between their pickup and delivery points and, therefore, rarely visit local hubs along their path. Since the (relative) difference between sorting and crossdocking time is smaller at access hubs, savings in handling time are also smaller (when sorting is avoided).

\begin{figure}[htp]
\begin{center}
  \includegraphics[width=\textwidth]{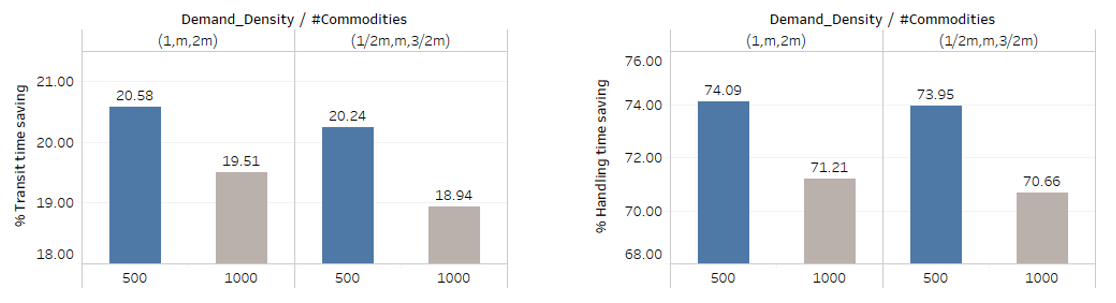}
  \caption{Sensitivity analysis on the number of commodities and the commodity sizes}
  \label{fig:numCom}
\end{center}
\end{figure}

Recall that to determine the size of a commodity, we use a triangular distribution with its mean, $m$, equal to the average volume per commodity (i.e., the total demand volume divided by the number of commodities). In Figure \ref{fig:numCom}, we show savings from containerization for different numbers of commodities and different commodity sizes induced by the minimum and maximum of the triangular distribution. 
A triangular distribution with parameters $(1,m,2m)$ implies that the minimum and maximum size of a commodity are 1 and $2m$, respectively. Similarly, a triangular distribution with parameters $(\frac{1}{2}m,m,\frac{3}{2}m)$ indicates that the minimum and maximum size of a commodity are 0.5 and 1.5 times the average volume, respectively.  To ensure a fair comparison, for a given total demand volume and number of commodities, we use the same set of commodities (i.e., the same origins, destinations, and service promises) but assign quantities from a triangular distribution with different parameters.

\begin{figure}[htp]
\begin{center}
  \includegraphics[width=\textwidth]{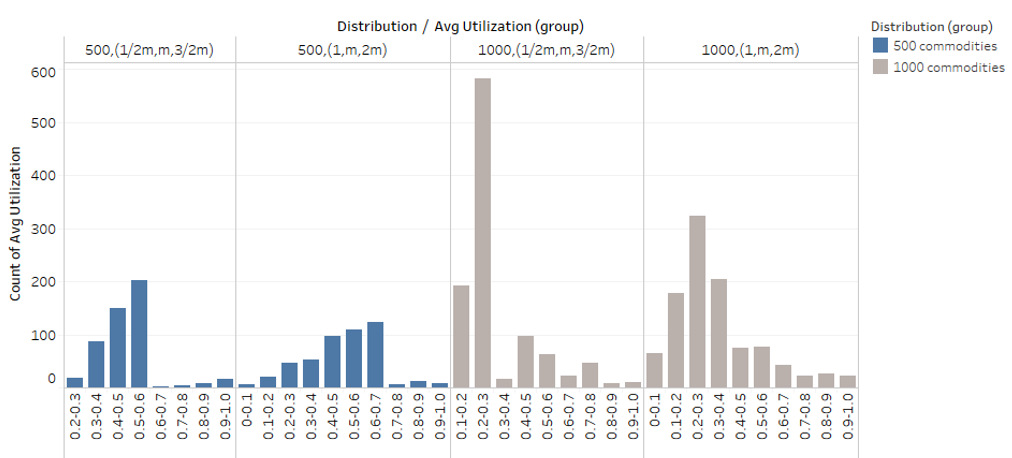}
  \caption{Container utilization for different demand size and distributions}
  \label{fig:distUtil}
\end{center}
\end{figure}

Comparing the results for the same demand volume but a different number of commodities (500 and 1000), we observe that when the average quantity per commodity is larger, the corresponding savings from containerization are more significant. This is intuitive, as more volume naturally travels together all the way from a commodity's origin to its destination. Moreover, a larger variation in commodity sizes implies higher handling time savings, translating into higher transit time savings from containerization. This impact is more pronounced when the commodity's average size is (relatively) small compared to the container size (e.g., when \# of commodities = 1000, demand volume = 10'000, and container size = 40). The relation between commodity size and handling time savings can be explained as follows. First, with a larger variation in commodity size, there is a higher chance to have commodities that almost fill a container by themselves and are transported from pickup point to delivery point in a container. Second, with a larger variation in commodity size, it is more likely that the remaining container space can be filled, resulting in higher container utilization (bin packing efficiencies). The container fill rate for the different number of commodities and commodity size distribution is illustrated in Figure \ref{fig:distUtil} which supports these arguments.

\subsection{Analysis on Different Network Configurations} 

In what follows, we examine the impact of the delivery network configuration on savings achieved through containerized consolidation. We evaluate the nine different network configurations shown in Figure \ref{fig:exprCityGrids} under different intracity demand patterns (i.e., uniform, centric, and bipolar). We create instances with 10,000 units of demand and 1000 commodities, with 50\% of commodities committed a 5-hour delivery service and the rest due in 10 hours. Regardless of the demand patterns, the intercity inbound, intercity outbound, and intracity commodities each represent a third of the total demand. We use the same set of commodities (with their origin, destination, size, and service) for all 9 network configurations to ensure a fair comparison. A commodity's alternative paths are generated considering up to 7 intermediate terminals and up to 5\% deviation from its shortest path. Lastly, we set the sorting to crossdocking time and capacity ratio at the hubs to 0.25 and 4, respectively. The IP model is solved to 0.01\% of optimality for all scenarios.

To facilitate the analysis of the results, we classify the nine network configurations using two characteristics. The first characteristic is the types of links between the hubs considered in each network configuration, called the ``Links Structure''. Based on this characteristic, the network structures fall into three categories: traditional hub and spoke (\textsc{HS}), hyperconnected 1 (\textsc{HC1}), and hyperconnected 2 (\textsc{HC2}). Recall that the difference between \textsc{HC1} and \textsc{HC2} is that in the former, each unit zone is linked to 4 access hubs, while in the latter, the access hubs tier is less dense, and each unit zone is served by only one access hub.
The second characteristic is the types of hubs included in each network configuration, called the ``Hubs Structure''. In the \textsc{Default} scenario, the intracity logistic network includes access, local, and gateway hubs, while in the \textsc{noLH} and \textsc{noGH} scenarios the local hub and gateway hub tiers are excluded, respectively. The different characteristics of the target network configurations are summarized in Table \ref{tab:networksChar}. 

\begin{table}[H]
   \scriptsize
   \centering
    \begin{tabular}{ll|ccc|ccc|c}
     \toprule
     \multicolumn{2}{c|}{\multirow{2}{*}{Network structure}} & \multicolumn{3}{c|}{Location (if present)} & \multicolumn{3}{c|}{Lateral shipment} & \multirow{2}{*}{\#AHs per zone}\\
    & & AH & LH & GH & AH & LH & GH & \\
    \midrule
    \multirow{3}{*}{\textsc{Default}} & HS & center & center & center & No & No & Yes & 1\\
    & HC1 & corner & corner & corner & Yes & Yes & Yes & 4\\
    & HC2 & corner & corner & corner & Yes & Yes & Yes & 1\\
    \midrule
    \multirow{3}{*}{\textsc{noLH}} & HS & center & - & center & No & - & Yes & 1\\
    & HC1 & corner & - & corner & Yes & - & Yes & 4\\
    & HC2 & corner & - & corner & Yes & - & Yes & 1\\
    \midrule
    \multirow{3}{*}{\textsc{noGH}} & HS & center & center & - & No & Yes & - & 1\\
    & HC1 & corner & corner & - & Yes & Yes  & - & 4\\
    & HC2 & corner & corner & - & Yes & Yes  & - & 1\\
    \bottomrule
    \end{tabular}
    \caption{Network Structures Characteristics} 
    \label{tab:networksChar}
\end{table} 

\begin{figure}[htp]
\begin{center}
  \includegraphics[width=0.8\textwidth]{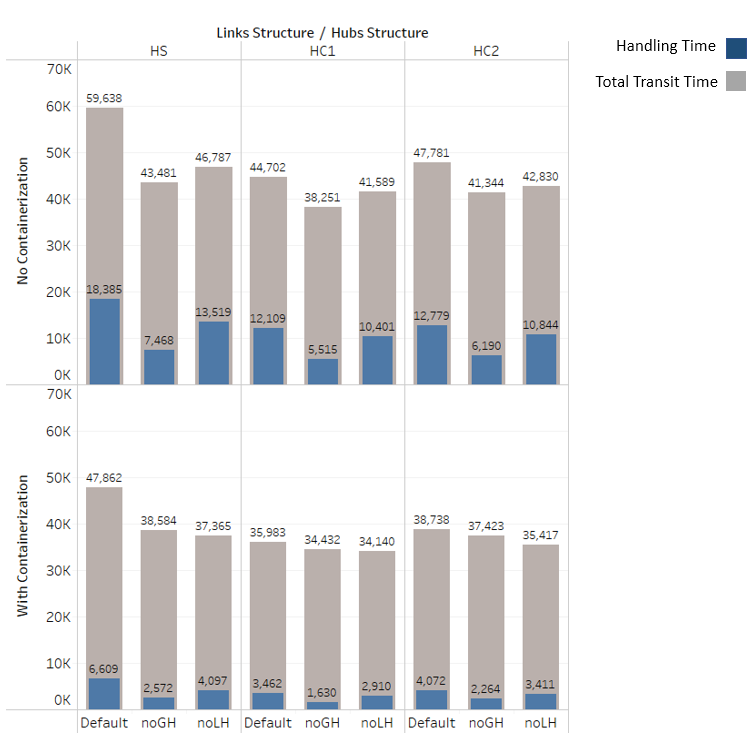}
  \caption{Total transit and handling time for uniform demand (COLORED)}
  \label{fig:uniform}
\end{center}
\end{figure}

Tables \ref{tab:uniDemand}, \ref{tab:centricDemand}, and \ref{tab:bipolarDemand} give the benefits of containerized consolidation (total transit and handling time savings) and the model solution times for different network configurations under uniform, centric, and bi-polar demand, respectively. The results indicate that regardless of the demand pattern and for any hub structure, \textsc{HS} results in the largest and \textsc{HC1} results in the smallest total transit and handling time both with and without containerization. Figure \ref{fig:uniform} shows these results for the uniform demand (for central and bi-polar demand refer to Appendix \ref{fig:centric} and \ref{fig:bipolar}). This happens because \textsc{HC1} provides the highest level of interconnection at the lower tiers of the network and, as such, does not force commodities to travel to the higher tiers unnecessarily.  Moreover, we see that for all three link structures, i.e., \textsc{HS}, \textsc{HC1} and \textsc{HC2}, the handling time and therefore the total transit time are the smallest when no gateway hubs are present.  This is intuitive because handling times at gateway hubs are larger than at the local and access hubs.

\begin{table}[htp]
   \scriptsize
   \centering
    \begin{tabular}{ll|rrr|rrr|rr}
     \toprule
     \multicolumn{10}{c}{Uniform Demand} \\
     \midrule
     \multicolumn{2}{c|}{\multirow{2}{*}{Network structure}} & \multicolumn{3}{c|}{total transit time} & \multicolumn{3}{c|}{handling time} & \multicolumn{2}{c}{solution time (sec)}\\
    & & noCont & withCont & \%Imprv & noCont & withCont & \%Imprv & noCont & withCont \\
    \midrule
    \multirow{3}{*}{\textsc{Default}} & HS & 59637 & 47862 & 19.74 & 18385 & 6609 & 64.05 & 0.70 & 237.97 \\ 
     & HC1 & 44702 & 35983 & 19.50 & 12109 & 3462 & 71.41 & 2.69 & 201.65 \\
     & HC2 & 47781 & 38738 & 18.92 & 12779 & 4072 & 68.14 & 1.02 & 22.38 \\ 
    \midrule
    \multirow{3}{*}{\textsc{noLH}} & HS & 46787 & 37365 & 20.14 & 13519 & 4097 & 69.70 & 0.63 & 4.48 \\ 
    & HC1 & 41589 & 34140 & 17.91 & 10401 & 2910 & 72.02 & 2.44 & 67.66 \\ 
    & HC2 & 42830 & 35418 & 17.31 & 10844 & 3411 & 68.54 & 0.81 & 5.83 \\ 
    \midrule
    \multirow{3}{*}{\textsc{noGH}} & HS & 43480 & 38584 & 11.26 & 7468 & 2572 & 65.56 & 1.25 & 28.56 \\ 
    & HC1 & 38251 & 34432 & 9.98 & 5515 & 1630 & 70.44 & 3.29 & 151.22 \\ 
    & HC2 & 41344 & 37423 & 9.48 & 6190 & 2264 & 63.43 & 0.97 & 19.01 \\ 
    \bottomrule
    \end{tabular}
    \caption{Containerized Consolidation savings across different network structures with Uniform Demand} 
    \label{tab:uniDemand}
\end{table} 

\begin{table}[H]
   \scriptsize
   \centering
    \begin{tabular}{ll|rrr|rrr|rr}
     \toprule
     \multicolumn{10}{c}{Centric Demand} \\
     \toprule
     \multicolumn{2}{c|}{\multirow{2}{*}{Network structure}} & \multicolumn{3}{c|}{total transit time} & \multicolumn{3}{c|}{handling time} & \multicolumn{2}{c}{solution time (sec)}\\
    & & noCont & withCont & \%Imprv & noCont & withCont & \%Imprv & noCont & withCont \\
    \midrule
    \multirow{3}{*}{\textsc{Default}} & HS & 54659 & 44140 & 19.25 & 16333 & 5814 & 64.41 & 0.61 & 82.20\\ 
    & HC1 & 42294 & 34192 & 19.16 & 11443 & 3446 & 69.88 & 2.45 & 90.03\\ 
    & HC2 & 44353 & 35889 & 19.08 & 11862 & 3758 & 68.32 & 0.91 & 15.32\\ 
    \midrule
    \multirow{3}{*}{\textsc{noLH}} & HS & 43204 & 35021 & 18.94 & 12078 & 3894 & 67.76 & 0.57 & 5.23 \\ 
    & HC1 & 39657 & 32921 & 16.99 & 9575 & 2816 & 70.59 & 2.43 & 61.01 \\ 
    & HC2 & 40632 & 33585 & 17.34 & 10287 & 3223 & 68.67 & 0.83 & 6.72 \\ 
    \midrule
    \multirow{3}{*}{\textsc{noGH}} & HS & 40273 & 35945 & 10.75 & 6750 & 2422 & 64.12 & 0.58 & 15.74 \\ 
    & HC1 & 36167 & 32791 & 9.34 & 5041 & 1632 & 67.63 & 2.78 & 85.85 \\ 
    & HC2 & 38073 & 34627 & 9.05 & 5444 & 1984 & 63.55 & 0.79 & 12.14 \\ 
    \bottomrule
    \end{tabular}
    \caption{Containerized Consolidation savings across different network structures with Centric Demand}
    \label{tab:centricDemand}
\end{table} 

\begin{table}[H]
   \scriptsize
   \centering
    \begin{tabular}{ll|rrr|rrr|rr}
     \toprule
     \multicolumn{10}{c}{Bi-polar Demand} \\
     \toprule
     \multicolumn{2}{c|}{\multirow{2}{*}{Network structure}} & \multicolumn{3}{c|}{total transit time} & \multicolumn{3}{c|}{handling time} & \multicolumn{2}{c}{solution time (sec)}\\
    & & noCont & withCont & \%Imprv & noCont & withCont & \%Imprv & noCont & withCont \\
    \midrule
    \multirow{3}{*}{\textsc{Default}} & HS & 59902 & 47477 & 20.74 & 18950 & 6524 & 65.57 & 0.73 & 1327.39 \\ 
    & HC1 & 46156 & 37339 & 19.10 & 12401 & 3760 & 69.68 & 4.30 & 543.50 \\ 
    & HC2 & 48566 & 39386 & 18.90 & 13234 & 4342 & 67.19 & 1.36 & 1025.94 \\ 
    \midrule
    \multirow{3}{*}{\textsc{noLH}} & HS & 46978 & 38099 & 18.90 & 13912 & 5034 & 63.82 & 0.67 & 20.32 \\ 
    & HC1 & 42235 & 34671 & 17.91 & 10745 & 3035 & 71.75 & 2.64 & 70.02 \\ 
    & HC2 & 42484 & 35451 & 16.55 & 10958 & 3884 & 64.56 & 0.88 & 15.65\\ 
    \midrule
    \multirow{3}{*}{\textsc{noGH}} & HS & 44356 & 39244 & 11.52 & 7969 & 2857 & 64.15 & 0.80 & 173.90\\ 
    & HC1 & 39861 & 35886 & 9.97 & 5843 & 1841 & 68.49 & 3.99 & 340.56\\ 
    & HC2 & 42183 & 38078 & 9.73 & 6676 & 2515 & 62.33 & 1.09 & 208.46\\ 
    \bottomrule
    \end{tabular}
    \caption{Containerized Consolidation savings across different network structures with Bi-polar Demand}
    \label{tab:bipolarDemand}
\end{table}

\begin{figure}[htp]
\centering
  \includegraphics[width=1.0\linewidth]{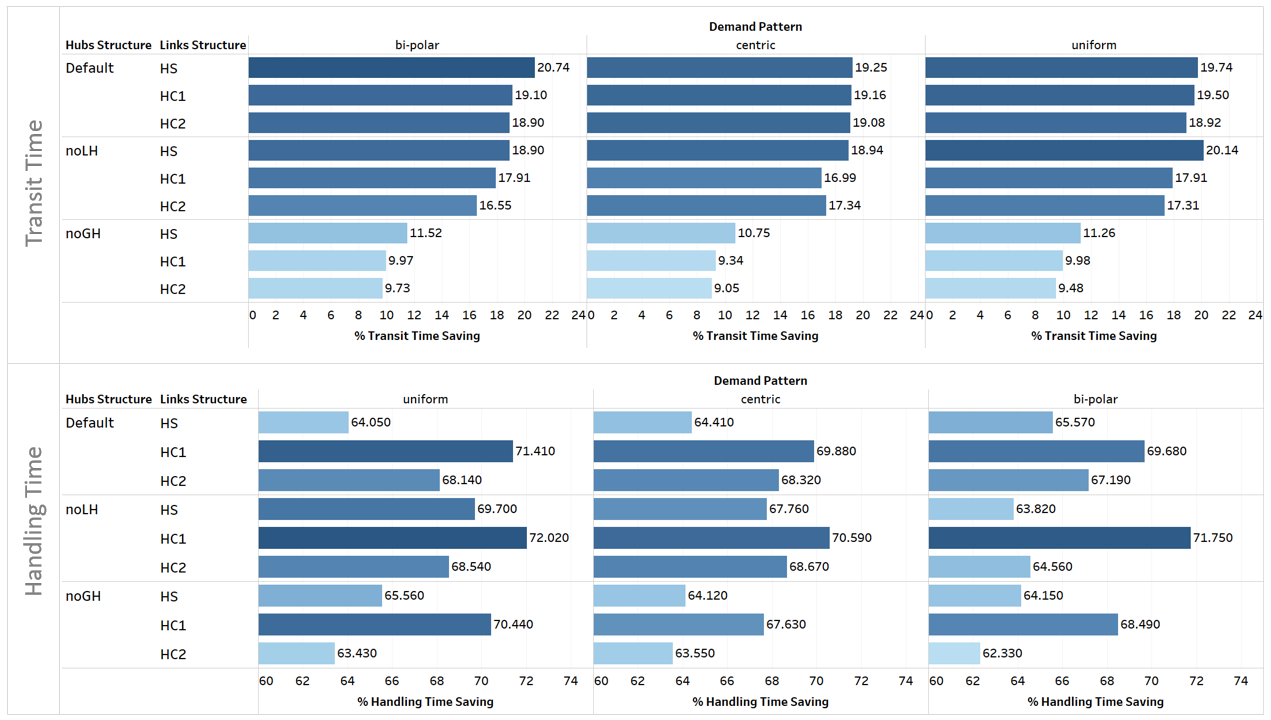}
  \caption{Total transit and handling time savings for uniform, centric, and bi-polar demand}
  \label{fig:timeSaving}
\end{figure}

Figures \ref{fig:timeSaving} show the total transit and handling time savings for the uniform, centric and bi-polar demand. As we see in this figure, when demand is uniform, regardless of the hub structure, the \textsc{HC1} links structure allows for the largest handling time savings because it provides more opportunities for commodities to merge and benefit from containerized consolidation. Moreover, we see that with a uniform demand, for the \textsc{Default} hub structure, the \textsc{HS} links structure results in the smallest handling time savings, while for the \textsc{noLH} and \textsc{noGH} hub structures, the \textsc{HC2} links structure results in the smallest handling time savings.

We observe a similar behavior when demand is centric, with one interesting distinction. With a uniform demand, when there is no local hub in the network, the total handling time savings is larger for the traditional \textsc{HS} structure in comparison to the \textsc{HC2} structure. That is while when demand is centric, chances for handling time savings are relatively higher in the \textsc{HC2} structure. The reason is, in the \textsc{HS} structure, each access hub is exclusively assigned to one unit zone, and therefore, when there are no local hubs in the network, this structure only offers one major point of flow consolidation (the single gateway hub) inside each urban area. Since with a centric demand, a large portion of commodities will travel within the same urban area, a smaller portion of total flow will bypass the gateway hubs sorting process, and as a result, relatively less handling time is saved. 

% transit time for uniform demand

Surprisingly, when demand is bi-polar, the \textsc{HS} and \textsc{HC2} structures with no local hubs allow for significantly less handling time savings compared to when demand is uniform. The reason is, when demand is bi-polar, these two structures make a large portion of commodities flow through few popular arcs in the middle of their paths. With such flow pooling, the model assigns less shipping capacity along the popular arcs and therefore enforces higher container utilization along those arcs. This may allow fewer commodities to get consolidated from \textit{near} their pickup to \textit{near} their delivery point and decrease the total handling time savings through containerized consolidation.

The percentage of transit time saving of a specific network configuration through containerized consolidation is impacted by the total handling time imposed along the commodities trip, as well as the potential for handling time savings. As such, a typical configuration may induce a relatively smaller handling time saving but larger transit time savings if the handling time constitutes a larger proportion of the commodities' total transit time. This phenomenon is specifically observed for the \textsc{HS} structure when compared to the other network configurations.

\begin{figure}[htp]
\centering
  \includegraphics[width=0.9\linewidth]{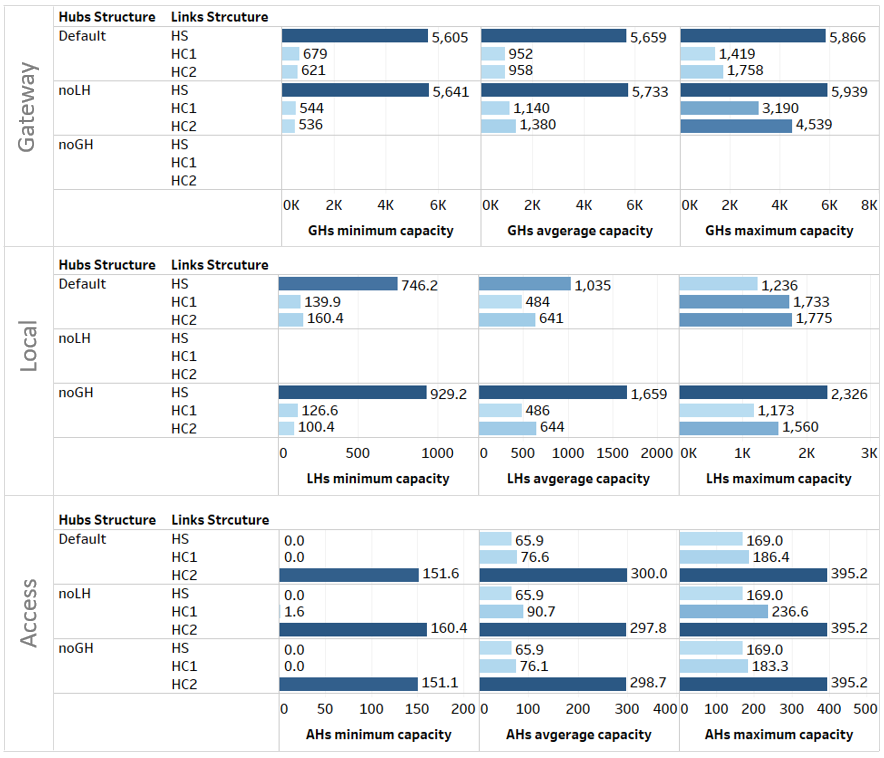}
  \caption{Gateway, local and access hubs capacity for different demand configurations}
  \label{fig:hubCap}
\end{figure}

Besides discussing the operational capabilities of different network configurations, that is, the total handling and transit times of commodities, it worth noting that each configuration may need different numbers and sizes of every hub type and, therefore, different levels of strategic investments. Figure \ref{fig:hubCap} shows the planned capacity based on the projected expected flow at each hub type for every network configuration when demand is uniform (see section \ref{capacity} for capacity planning methodology).

As we see in Figure \ref{fig:hubCap}, for every hub structure, the minimum, average, and maximum capacity for the access hubs are significantly larger in the \textsc{HC2} link structure, when compared to \textsc{HS} and \textsc{HC1}. This happens because \textsc{HC2} has smaller number of access hubs (64 as opposed to 256 and 289 for \textsc{HS} and \textsc{HC1} structures, respectively), and every 4 unit zones are served by only one access hub. The \textsc{HC1} link structure stands in second place in terms of the maximum access hubs capacity, but its difference with the \textsc{HS} structures is not significant.

For the local hubs, the minimum and average capacities are the largest in the \textsc{HS} structure. This happens because most of the commodities are forced to travel to the higher tiers of the network and, therefore, will pass through the local hubs. Interestingly, in terms of maximum local hub capacity, \textsc{HS} stands in the first place when no gateway hubs exist and the last place when all hub types are present. The reason is that when gateway hubs are present, the local hubs located next to the gateway hubs in \textsc{HC1} and \textsc{HC2} structures become very popular as they can bridge the flow to the highest tier of the network. For the same reason, when no gateway hubs are present, the maximum flow passing through the local hubs significantly decreases for the \textsc{HC1} and \textsc{HC2} link structures, while it increases for the \textsc{HS} structure.

Lastly, in terms of gateway hubs' capacity, the \textsc{HS} link structure requires a significantly larger minimum, average, and maximum capacity at the gateway hubs. This again happens since in this structure, most of the flow is directed to the highest tier of the network and passes through the gateway hubs. \textsc{HC1} requires the least average and maximum capacity at the gateway hubs as it allows for the highest level of interconnection between access and local hubs in the lower tiers of the network. Therefore, when no local hubs are present, the average and maximum capacity required at the gateway hubs increases in the \textsc{HC1} and \textsc{HC2} structures, yet staying far below the capacity required at the \textsc{HS} link structure.

%section 6
\section{Concluding remarks} \label{sec:conclusion}

We have demonstrated that significant benefits can be achieved, in terms of total in-transit and handling time of commodities, through joint parcel routing and containerized consolidation in megacity parcel logistics. We have also extensively investigated the impact of different strategic, tactical, and operational characteristics of a logistic system on the potential benefits of containerized consolidation; Such characteristics include the logistic network configuration, demand pattern, and sorting capacity and capability, to name a few.

As this is the first study of its kind and for the sake of simplicity, we have considered a single container size. The natural next step is to explore the benefits of multiple container sizes to the crossdocking rate, containers' utilization, and the total shipping volume requirements. Even though extending the proposed integer programming model to accommodate multiple container sizes is fairly straightforward, the solution of instances of meaningful size will be significantly more difficult and more sophisticated, and customized solution approaches will have to be developed. This is left for future research. 

Another possible research avenue to explore in the future is the trade-off between cost and service. Providing additional capacity (hub and link capacity) will increase cost but will also likely improve service. A better understanding of this trade-off curve will be valuable for companies operating in the last-mile logistics space. 

Lastly, this study has approached the containerized consolidation model as a tactical problem. However, developing smart heuristic, meta-heuristic or intelligent hybrid approaches are also encouraged to better deal with the dynamic uncertainty of urban parcel logistics and its higher pressure on the solution time.
  %section 7

%######### References ##########
\newpage
\bibliography{refs.bib} 
\bibliographystyle{ieeetr}

%########## Appendix #########
%\newpage
\newpage
\section*{Appendices}
\appendix

\counterwithin{figure}{section}

%####################################
%####################################
%####################################
\section{Total transit and handling time}\label{app_A}

% In-transit time for centric demand
\begin{figure}[H]
\begin{center}
  \includegraphics[width=0.85\textwidth]{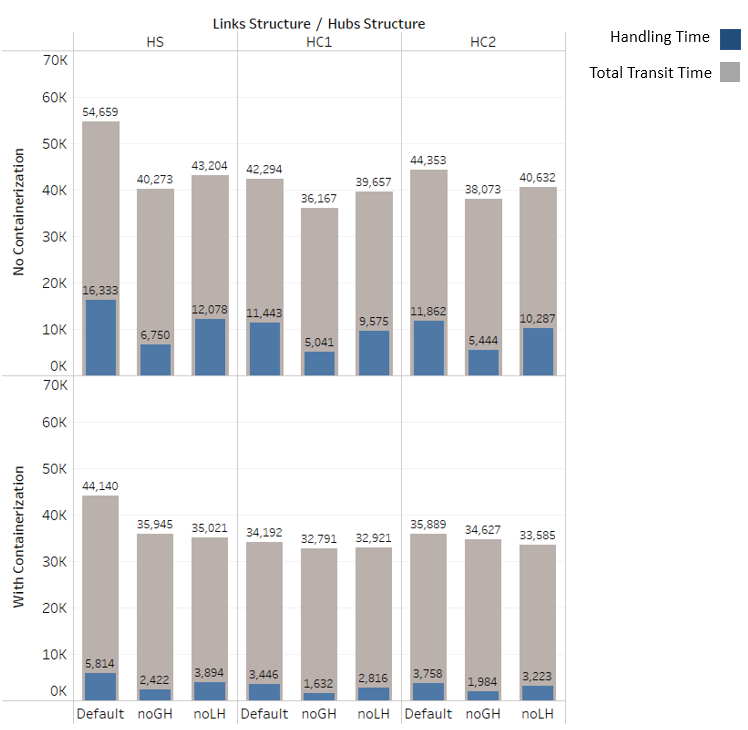}
  \caption{Total transit and handling time for centric demand (COLORED)}
  \label{fig:centric}
\end{center}
\end{figure}

% In-transit time for bi-polar demand
\begin{figure}[H]
\begin{center}
  \includegraphics[width=0.85\textwidth]{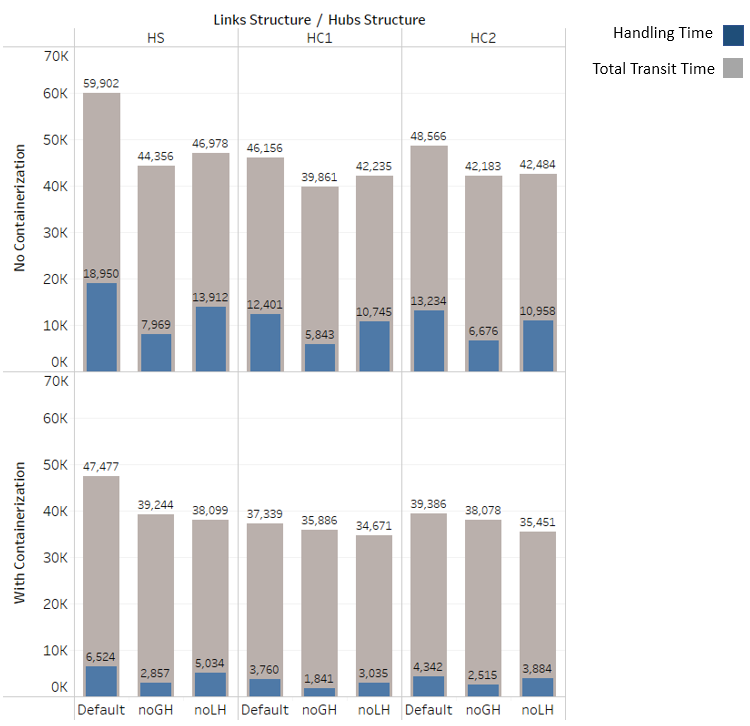}
  \caption{Total transit and handling time for bi-polar demand (COLORED)}
  \label{fig:bipolar}
\end{center}
\end{figure}
\end{document}